\documentclass[12pt]{amsart}
\newcommand{\p}{\mathbb P}

\newcommand{\lra}{\longrightarrow}
\newcommand{\mo}{\mathcal O}
\newcommand{\tp}{\tilde\mathbb P^3}
\newcommand{\ml}{\mathcal L}

\newcommand{\vd}{$\vdots$}
\theoremstyle{plain}
\newtheorem{thm}{Theorem}[section]
\newtheorem{prop}[thm]{Proposition}
\theoremstyle{definition}
\newtheorem{defn}{Definition}[section]
\author{S\l awomir Cynk}
\thanks{Partially supported by KBN grant no 2P03A 083 10}
\subjclass{\mbox{Primary 14J32,14J17;
Secondary 14B05}}
\keywords{Calabi--Yau manifolds, double solids, surface singularities}
\title[Double  octic arrangements with isolated
singularities]
{Double coverings of octic arrangements with isolated singularities}
\address{S\l awomir Cynk\\ Institute of Mathematics\\ Jagiellonian University\\
  ul. Reymonta 4\\ 30-059 Krak\'ow\\ Poland} 
\email{cynk@im.uj.edu.pl}
\begin{document}
\begin{abstract}
In this paper we construct 206 examples of Calabi--Yau manifolds with
different Euler numbers. All constructed examples are smooth models of double
coverings of $\p^{3}$ branched along an octic surface. 
We allow 11 types of (not necessary isolated)
singularities in the branch locus. Thus we
broaden the class of examples studied in \cite{CS}.
For every considered example we compute the Euler number and give a precise
description of a resolution of singularities.
\end{abstract}
\maketitle
\section{Introduction}

\label{s:intr}
In this paper we study a class of Calabi--Yau manifolds. By a Calabi--Yau
manifold we mean a k\"ahler, smooth threefold with trivial canonical bundle
and no global 1-forms. One of the method of constructing Calabi--Yau 
manifold is to study a double covering of $\p^{3}$ branched along an 
octic surface.

Let $B$ be a surface in $\p^3$ of degree eight. Let
$\ml=\mo_{\p^3}(4)$, $\ml$ is a line bundle on $\p^3$ such that
$\ml^{\otimes 2}=\mo_{\p^3}(B)$. Then $\ml $ defines a double covering of
$\p^3$ 
branched along $B$ and the singularities of the double cover are in
one--to--one correspondence with singularities of $B$. 

If $B$ is smooth then the resulting double covering is a Calabi--Yau
threefold, if $B$ has only nodes (ordinary double points) then the
double cover has also only nodes,  and these nodes can be resolved by
the mean of small resolution. In this case again the resulting
threefold is Calabi--Yau. This construction was precisely described in
\cite{Clemens}.

In \cite{CS} a case of $B$ being an octic arrangement (i.e. a surface
which locally looks like a plane arrangement) was studied.
In this paper we shall use methods introduced in~\cite{CS} 
to study a bigger class of octic surfaces,
namely we shall consider arrangements with ordinary multiple points of
multiplicity $2$, $4$ and $5$. 
Altogether we allow 11 types of singularities of branch locus.

Our main result is the following theorem
\begin{thm}

\label{t:main}
\samepage
If an octic arrangement $B$ contains only 
\begin{itemize}
\item double and triple curves,
\item arrangement $q$--fold points, $q=\,2,\,3,\,4,\,5$,
\item isolated $q$--fold points $q=\,2,\,4,\,5$
\end{itemize}
then the double covering of $\p^3$ branched along $B$ has a
non--singular model $\hat X$ which is a Calabi--Yau
threefold. 

Moreover if $B$ contains no triple elliptic curves then
\begin{multline*}
\mathbf e(\hat
X)=8-\sum_i(d_i^3-4d_i^2+6d_i)+2\sum_{i<j}(4-d_i-d_j)d_id_j-
\sum_{i<j<k}d_id_jd_k\\
+4p^0_4+3p^1_4+16p^0_5+18p^1_5+20p^2_5+l_3+2m_2+36m_4+56m_5.
\end{multline*}
\end{thm}
The idea of the proof of this theorem is to give a resolution of singularities
of $X$ by a sequence of admissible blowing-ups (i.e. blowing-ups that do not
affect the first Betti number and the canonical divisor of the double
covering, cf.~\cite{CS}).

We apply Theorem~\ref{t:main} to give examples of 
Calabi--Yau manifolds with 206
different Euler number (we realize any even number from the interval
$\left<-296,104\right>$ as an Euler number of a Calabi--Yau manifold).

\section{Surfaces with ordinary multiple points}

\label{s:top}
Let $B$ be a surface in $\p^3$ with only ordinary multiple points. 
That means that if we 
 consider the blowing--up $\sigma:\tilde\p^3\lra\p^3$ of $\p^3$
at all singular points of $B$ then  the strict
transform $\tilde B$ of $B$  is smooth and intersect the exceptional
divisor of $\sigma$ transversally. 
Let $m_p$ denotes the number of $p$--fold points 
on $B$. The following Proposition contains the numerical
data of $B$ and $\tilde B$.
\begin{prop}
\label{p:surf}
\begin{eqnarray}
c_1^2(\tilde B)&=&d(d-4)^2-\sum_p(p-2)^2pm_p\\
c_2(\tilde B)&=&d^3-4d^2+6d-\sum _p(p-2)p^2m_p\\
e(B)&=&d^3-4d^2+6d-\sum _p(p-1)^3m_p\\
p_a(\tilde B)&=&\binom{d-1}{3}-\sum_p\binom p3m_p
\end{eqnarray}
\end{prop} 

\section{Octic arrangements with isolated singularities}

\begin{defn}\samepage
An {\em octic arrangement with isolated
singularities}  is a
surface $B\subset \p^3$ of degree 8 which is a sum of irreducible
surfaces $B_1,\dots,B_r$ with only ordinary multiple points which 
satisfies the following conditions:

\begin{enumerate}
\item For any $i\not=j$ the surfaces $B_i$ and $B_j$ intersects
transversally along a smooth irreducible curve $C_{i,j}$ or they are
disjoint, 
\item The curves $C_{i,j}$ and $C_{k,l}$ are either equal or disjoint or they
  intersects transversally.  
\end{enumerate}
\end{defn}
This definition is a generalization of the notion of octic arrangement
introduced in \cite{CS} (where the surfaces $S_i$ are assumed to  be
smooth). Observe that from (1) Sing$B_i\cap B_j=\emptyset$ for
$i\not=j$. We shall denote $d_i:=\dim B_i$.

A singular point of $B_i$ we shall call an {\em isolated singular point  of
the arrangement}.
A point $P\in B$ which belongs to $p$  of surfaces $B_1,\dots,B_r$ we
shall call an {\em arrangement $p$--fold point}.
We say that irreducible curve $C\subset B$ is a
$q$--fold curve if exactly $q$ of surfaces $B_1,\dots,B_r$ passes
through it. 

We shall use the following numerical data for an arrangement:
\begin{itemize}
\item[$p^i_q$] Number of arrangement $q$--fold points lying on exactly
$i$ triple curves,
\item[$l_3$] Number of triple lines, 
\item[$m_q$] number of isolated $q$--fold points.
\end{itemize}

\section{Proof of Theorem~\ref{t:main}}
\label{sec:proof}
Let $B$ be an octic arrangement satisfying assumptions of the Theorem and let
$X$ be a double covering of $\p^{3}$ branched along $B$. We shall
find a sequence of admissible blowing--ups (i.e. blowing--ups of double or
triple curves and 4--fold or 5--fold points)
$\sigma:\p^{\ast}\longrightarrow\p^{3}$  
and a reduced divisor $B^{\ast}\subset \p^{\ast}$ with ordinary double points (nodes) as the
only singularities and such that 
\begin{itemize}
\item []$\tilde B\le B\le B^{\ast}$ (where $\tilde B$ is a strict transform
  and $B^{\ast}$ is a pullback of $B$ by $\sigma$),
\item []$B^{\ast}$ is even as an element of the Picard group Pic$\;\p^{\ast}$.
\end{itemize}

Let us now describe an algorithm to obtain $\sigma$,
the method is in fact a  modification of the method
introduced in \cite{CS}. We resolve all
singularities of $B$ except the nodes

{\bf 1. Resolution of multiple curves and arrangement multiple
points.} In these cases we shall apply the method described in
\cite{CS}. 

{\bf 2. Resolution of isolated $4$--fold points}. 
We  blow--up a $4$--fold point, and then replace the branch
locus by its strict transform.

{\bf 3. Resolution of isolated $5$--fold points}. We  blow--up a
$5$--fold point,  and then replace the branch
locus by its strict transform plus the exceptional divisor. The proper
transform intersects exceptional divisor transversally along a smooth
plane curve of degree 5. We treat this curve in the same way as as 
an arrangement double curve
i.e. we blow-up this curve, and replace the branch locus by its strict
transform.

The double covering of $\p^{\ast}$ branched along $B^{\ast}$ has nodes
(corresponding to nodes of $B^{\ast}$) as the only singularities. 

{\bf 4. Resolution of nodes.}
There are two possibilities for the resolution of a node
on a 3--dimensional variety: blow--up or small-resolution (for details
see \cite{Clemens})

We shall denote the non--singular model of $X$ by $\tilde X$ if in step 4 we
choose a blow--up and by $\hat X$ if we choose a small resolution.
To any node on $B$  we have associated a line on $\hat X$. $\tilde X$
is a blowing-up of  $\hat X$ at  all this lines. As a
consequence we see that
$e(\tilde X)=e(\hat X)+2m_2$.

The blow-ups  used in steps 1---3 are (according to \cite{CS})
admissible, i.e. they do not affect the first Betti number and
canonical divisor of $\tilde X$. We see therefore that 
\[K_{\tilde X}=E_2\]
where $E_2$ denotes the exceptional divisor on $\tilde X$ associated
to all nodes of $B$ hence \[K_{\hat  X}=0.\]

In order two compute $e(\tilde X)$ we compare our case with the one
studied in \cite{CS}. From the Proposition~\ref{p:surf} we see that in
our case $e(\tp)$ increases by $2m_2+2m_4-8m_5$ whereas
$e(\tilde B)$ decreases by $32m_4+72m_5$.

The Euler number $e(\hat X)$ is hencefore greater by 
\[2m_2+36m_4+56m_5\]
in comparison with the case with no isolated singular points.
Using \cite[Thm.~3.4]{CS} proves the theorem.
\qed

\def\arraystretch{1.13}
\begin{center}
\begin{table}[p]
\begin{tabular}{||l||r|r|r|r|r|r|r||r||}
\hline\hline
$(d_1,d_2,\dots,d_r)$&$p_4^0$&$p_4^1$&$p_5^0$&$p_5^1$&$
p_5^2$&$l_3$&$m_2$&$\mathbf e(\hat
X)$
\\\hline\hline
8                               &       &       &       &       &       &       &       &-296\\\cline{2-9}
8                               &       &       &       &       &       &       &1      &-294\\\cline{2-9}
\vd                             &\vd    &       &       &       &       &       &\vd    &\vd    \\\cline{2-9}     
8                               &       &       &       &       &       &       &107    &-82    \\\hline  
(1,1,2,4)                       &1      &       &       &       &       &       &       &-80    \\\cline{2-9}
(1,1,2,4)                       &1      &       &       &       &       &       &1      &-78    \\\cline{2-9}
\vd                             &\vd    &       &       &       &       &       &\vd    &\vd    \\\cline{2-9}
(1,1,2,4)                       &1      &       &       &       &       &       &16     &-48    \\\hline
(1,1,1,1,4)             &       &       &1      &       &       &       &1      &-46    \\\cline{2-9}
\vd                             &\vd    &       &       &       &       &       &\vd    &\vd    \\\cline{2-9}
(1,1,1,1,4)             &       &       &1      &       &       &       &16     &-16    \\\hline
(1,1,1,1,1,3)           &       &       &       &       &       &       &1      &-14    \\\cline{2-9}
\vd                             &\vd    &       &       &       &       &       &\vd    &\vd    \\\cline{2-9}
(1,1,1,1,1,3)           &       &       &       &       &       &       &4      &-8     \\\cline{2-9}   
(1,1,1,1,1,3)           &2      &       &       &       &       &       &1      &-6     \\\cline{2-9}
\vd                             &\vd    &       &       &       &       &       &\vd    &\vd    \\\cline{2-9}
(1,1,1,1,1,3)           &2      &       &       &       &       &       &4      &0      \\\cline{2-9}
(1,1,1,1,1,3)           &       &5      &       &       &       &1      &1      &2      \\\cline{2-9}
\vd                             &\vd    &       &       &       &       &       &\vd    &\vd    \\\cline{2-9}
(1,1,1,1,1,3)           &       &5      &       &       &       &1      &4      &8      \\\hline
(1,1,1,1,2,2)           &2      &       &       &       &       &       &1      &10     \\\hline
(1,1,1,1,1,3)           &       &3      &       &1      &       &1      &       &12     \\\cline{2-9}
\vd                             &\vd    &       &       &       &       &       &\vd    &\vd    \\\cline{2-9}
(1,1,1,1,1,3)           &       &3      &       &1      &       &1      &4      &20     \\\hline
(1,1,1,1,1,1,2) &       &       &       &       &       &       &1      &22     \\\hline
(1,1,1,1,1,3)           &       &1      &       &2      &       &1      &       &24     \\\cline{2-9}
\vd                             &\vd    &       &       &       &       &       &\vd    &\vd    \\\cline{2-9}
(1,1,1,1,1,3)           &       &1      &       &2      &       &1      &4      &32     \\\hline
(1,1,1,1,1,1,2) &3      &       &       &       &       &       &1      &34     \\\cline{2-9}
(1,1,1,1,1,1,2) &       &       &1      &       &       &       &       &36     \\\cline{2-9}
(1,1,1,1,1,1,2) &       &       &1      &       &       &       &1      &38     \\\hline
(1,1,1,1,1,1,1,1)       &       &       &       &       &       &       &       &40     \\\hline
(1,1,1,1,1,1,2) &1      &5      &       &       &       &1      &1      &42
\\\hline\hline 

\end{tabular}
\end{table}
\begin{table}[p]
\begin{tabular}{||l||r|r|r|r|r|r|r||r||}
\hline\hline
$(d_1,d_2,\dots,d_r)$&$p_4^0$&$p_4^1$&$p_5^0$&$p_5^1$&$
p_5^2$&$l_3$&$m_2$&$\mathbf e(\hat
X)$\\\hline\hline
(1,1,1,1,1,1,1,1)       &1      &       &       &       &       &       &       &44     \\\hline
(1,1,1,1,2,2)           &1      &1      &       &2      &       &1      &1      &46     \\\hline
(1,1,1,1,1,1,1,1)       &2      &       &       &       &       &       &       &48     \\\hline
(1,1,1,1,1,1,2) &       &3      &       &1      &       &1      &1      &50     \\\hline
(1,1,1,1,1,1,1,1)       &3      &       &       &       &       &       &       &52     \\\hline
(1,1,1,1,1,1,2) &1      &3      &       &1      &       &1      &1      &54     \\\hline
(1,1,1,1,1,1,1,1)       &4      &       &       &       &       &       &       &56     \\\hline
(1,1,1,1,1,1,2) &2      &3      &       &1      &       &1      &1      &58     \\\hline
(1,1,1,1,1,1,1,1)       &5      &       &       &       &       &       &       &60     \\\hline
(1,1,1,1,1,1,2) &       &6      &       &       &1      &2      &1      &62     \\\hline
(1,1,1,1,1,1,1,1)       &6      &       &       &       &       &       &       &64     \\\hline
(1,1,1,1,1,1,2) &1      &6      &       &       &1      &2      &1      &66     \\\hline
(1,1,1,1,1,1,1,1)       &7      &       &       &       &       &       &       &68     \\\hline
(1,1,1,1,1,1,2) &2      &6      &       &       &1      &2      &1      &70     \\\hline
(1,1,1,1,1,1,1,1)       &8      &       &       &       &       &       &       &72     \\\hline
(1,1,1,1,1,1,2) &       &4      &       &1      &1      &2      &1      &74     \\\hline
(1,1,1,1,1,1,1,1)       &9      &       &       &       &       &       &       &76     \\\hline
(1,1,1,1,1,1,2) &1      &4      &       &1      &1      &2      &1      &78     \\\hline
(1,1,1,1,1,1,1,1)       &       &1      &       &2      &       &1      &       &80     \\\hline
(1,1,1,1,1,1,2) &2      &4      &       &1      &1      &2      &1      &82     \\\hline
(1,1,1,1,1,1,1,1)       &       &8      &       &1      &       &2      &       &84     \\\hline
(1,1,1,1,1,1,2) &       &2      &       &2      &1      &2      &1      &86     \\\hline
(1,1,1,1,1,1,1,1)       &12     &       &       &       &       &       &       &88     \\\hline
(1,1,1,1,1,1,2) &1      &2      &       &2      &1      &2      &1      &90     \\\hline
(1,1,1,1,1,1,1,1)       &       &4      &       &1      &1      &2      &       &92     \\\hline
(1,1,1,1,1,1,2) &2      &2      &       &2      &1      &2      &1      &94     \\\hline
(1,1,1,1,1,1,1,1)       &       &6      &       &2      &       &2      &       &96     \\\hline
(1,1,1,1,1,1,2)         &3      &2      &       &2      &1      &2      &1      &98     \\\cline{2-9}
(1,1,1,1,1,1,2)       &4      &2      &       &2      &1      &2      &       &100    \\\cline{2-9}
(1,1,1,1,1,1,2)         &4      &2      &       &2      &1      &2      &1      &102    \\\hline
(1,1,1,1,1,1,1,1)       &       &7      &       &       &2      &3      &       &104    \\\cline{2-9}
(1,1,1,1,1,1,1,1)       &       &9      &       &1      &1      &3      &       &108    \\\cline{2-9}
(1,1,1,1,1,1,1,1)       &       &3      &       &       &3      &3      &       &112    \\\cline{2-9}
(1,1,1,1,1,1,1,1)       &1      &3      &       &       &3      &3      &       &116    \\\cline{2-9}
(1,1,1,1,1,1,1,1)       &2      &3      &       &       &3      &3      &       &120    \\\cline{2-9}
(1,1,1,1,1,1,1,1)       &       &4      &       &       &4      &4      &       &136    \\\hline\hline

\end{tabular}
\end{table}
\end{center}
\section{Examples}
\def\arraystretch{1}
In this section we shall apply Theorem\ref{t:main} to study various examples
of octic arrangement. As a result we obtain 206 examples of Calabi--Yau
manifolds with different Euler
numbers, we shall for instance realize every even number from the interval
$\left<-296,104\right>$ as an Euler number of a Calabi--Yau manifold.

We shall need information about the number of nodes allowed on a nodal surface
of degree $\le 8$ in $\mathbb P^{3}$
Using results from \cite{Ba, B, Catanese}
we can formulate the following proposition
\begin{prop}\rule{0pt}{0pt}
\label{p:nod}\samepage
\begin{itemize}

\item[a.] For $m_2=0, 1, 2, 3, 4$ there
exists a nodal cubic surface with
exactly $m_2$ nodes,
\item[b.] For $m_2=0, 1, \dots, 16$ there
exists a nodal quartic surface with
exactly $m_2$ nodes,
\item [c.] For $m_2=0, 1, \dots, 65$ there
exists a nodal sextic surface with
exactly $m_2$ nodes,
\item [d.] For $m_{2}=0,1,\dots,107$ there 
exists a nodal octic surface with 
exactly $m_{2}$ nodes. 
\end{itemize}
\end{prop}

Using the above Proposition and Theorem 1.1 we can compile a table containing
numerical data of octic arrangements and corresponding
Euler numbers. 
Most of Euler numbers can be obtained from several different
arrangements, in the table we give one example per number. 
In the table we avoid arrangements with 4--fold and 5--fold points, they do not
leave to new Euler numbers. On the other hand arrangements with 4--fold and
5--fold points usually have higher Picard number then the ones with only
nodes.

Most of the examples in the table are modification of the ones given in
\cite{CS} obtained by adding isolated singularities.
In many cases it is easy to write down explicit equation of the branch
locus. Proof of Theorem~\ref{t:main} gives a detailed description of a
resolution of singularities. Although the resolution of singularities is not
uniquely determined, the different resolutions of the same double solid differ
only by flop. Consequently most of the numerical data (like f.i. Euler number)
are uniquely determined (cf.~\cite{kollar}).

\subsection*{Acknowledgment}\parindent=0cm
Part of this work was done during the Author stay at the Erlangen--N\"unmberg
University supporetd by DFG (project number 436 POL 113/89/0). I would
like to thank Prof. W. Barth for suggesting the problem and 
valuable remarks.

\end{document}